\documentclass[letterpaper, 10 pt, conference,compsoc]{ieeeconf}  

\IEEEoverridecommandlockouts                              
\usepackage[T1]{fontenc}
\usepackage[utf8]{inputenc}
\usepackage{comment}
\usepackage{palatino} 
\usepackage{mathrsfs}
\usepackage{graphicx}
\usepackage{pstool} 
\usepackage{mathtools}
\usepackage{amsmath} 
\usepackage{amssymb}  
\usepackage{color}
\usepackage{enumerate}
\usepackage{multirow}
\usepackage{nomencl}
\usepackage{array,url}
\usepackage{bm}
\usepackage{bigints}
\usepackage[acronym]{glossaries}
\newtheorem{theorem}{\bf Theorem}[section]

\newtheorem{proposition}[theorem]{\bf Proposition}

\newtheorem{remark}{\bf Remark}[section]

\newtheorem{problem}{\bf Problem}

\newcommand{\bmat}[1]{\begin{bmatrix}#1\end{bmatrix}}
\newacronym{hjbe}{HJBE}{Hamilton-Jacobi-Bellman equation}
\newacronym{nhim}{NHIM}{normally hyperbolic invariant manifold}
\title{\LARGE \bf
Maximally Degenerate Floquet Structure and Possible Nonexistence of Optimal Control in a Pendulum Swing-Up Problem
}

\author{Noboru Sakamoto${}^{1,\ast}$ 
\thanks{1: Faculty of Science and Technology, Nanzan University, Yamazato-cho 18, Showa-ku, Nagoya, 464-8673, Japan}%
\thanks{$\ast$:  
{\tt\footnotesize E-mail:  noboru.sakamoto@nanzan-u.ac.jp}
}%
}%
\makeatletter
\providecommand\HyPL@Entry[1]{}
\AddToHook{env/document/begin}{%
\@ifpackageloaded{preview}{
\ifPreview
 \let\Hy@FirstPageHook\relax
 \let\Hy@EveryPageAnchor\relax
\fi}{}}
\makeatother
\begin{document}
\maketitle
\thispagestyle{empty}
\pagestyle{empty}

\begin{abstract}
This paper studies an infinite-horizon optimal control problem for a pendulum with quadratic control cost via its associated Hamiltonian system. The problem is strongly degenerate, as the linearization at the upright equilibrium has purely imaginary eigenvalues with multiplicity, making standard linear analysis inconclusive.
Using Lie series and Chetaev’s instability theorem, we show that the equilibrium is weakly unstable in a non-hyperbolic sense. We further identify a family of periodic orbits forming an invariant cylinder, whose monodromy matrix exhibits a maximally degenerate Floquet structure consisting of a single Jordan block at the unit multiplier.
This structure implies slow, non-exponential dynamics and provides a mechanism for long-time transitions with arbitrarily small control effort. Although derived for a pendulum, this phenomenon arises from degeneracy in the optimal control formulation and may occur more broadly, suggesting possible nonexistence of optimal control despite stabilizability.
\end{abstract}
\section{Introduction} 
In control science, the pendulum is a benchmark that serves as an educational example and as a testbed for advanced control design, including stabilization \cite{Angeli:01:automatica}, under-actuated systems, and swing-up control \cite{Astrom:00:automatica,Shiriaev:01:syscon,Astrom:08:automatica,LaHera:09:mechatronics}. 
Controlling the nonlinear motion of a pendulum with optimal control is also studied by many authors, and in \cite{Horibe:17:ieee_cst}, an experimental result of optimal swing-up is reported, approximately solving a \gls{hjbe}.

Although the existence of optimal control can, in principle, be characterized via solutions of the \gls{hjbe}, this approach is not suitable for determining a priori whether an optimal control exists for a given nonlinear system, even though it is effective for constructing approximate optimal controls, since analytic solutions are generally unavailable. In \cite{Sakamoto:22:automatica}, the author established a sufficient condition for the existence of an infinite-horizon optimal control for a class of nonlinear systems based on stabilizability, detectability, and growth conditions. 

On the other hand, in \cite{Horibe:17:ieee_cst}, it is reported that there exist multiple (locally) optimal controls for the problem of pendulum swing-up. The aforementioned paper \cite{Sakamoto:22:automatica} applies its main result to show that an optimal control exists for this problem if the cost functional satisfies a detectability condition. If the cost is not detectable, however, a highly nonstandard behavior may arise, as is illustrated in \S~\ref{sctn:simulation} of the present paper, showing the possibility of an infinite number of local optimal controls. 

In this paper, we analyze the Hamiltonian system associated with the infinite-horizon optimal control problem for the pendulum in the absence of detectability, and discover a highly nonstandard dynamical structure. 
A key step in this analysis is to understand the stability properties of the equilibrium of an associated Hamiltonian system corresponding to the hanging position. Instability of this equilibrium is essential for the existence of trajectories departing from it and forming connecting orbits. However, due to the strong degeneracy of the linearized system, standard linearization-based methods are inconclusive. Therefore, we employ nonlinear tools, including Lie series and Chetaev’s instability theorem, to characterize the local dynamics. As will be shown, the instability is weak and non-hyperbolic, reflecting an underlying degeneracy that is further revealed by the monodromy analysis. 
We further show that a family of periodic orbits forms an invariant cylinder whose monodromy matrix exhibits a maximally degenerate Floquet structure, consisting of a single Jordan block associated with the unit multiplier. This degeneracy leads to a weak (non-hyperbolic) instability mechanism and suggests the possibility of long-time trajectories achieving large state transitions with arbitrarily small control effort. This observation is closely related to the existence theory of infinite-horizon optimal control developed in \cite{Sakamoto:22:automatica}, and indicates that, in the present problem, optimal controls may fail to exist despite stabilizability.

Although the analysis is carried out for a pendulum, the phenomenon revealed here is not specific to this system. The degeneracy arising from the absence of state weighting (Q = 0) leads to non-hyperbolic dynamics in the associated Hamiltonian system, which may occur in a broader class of optimal control problems, such as nonlinear mechanical systems. The present result provides a concrete example where such degeneracy leads to slow dynamics and potential nonexistence of optimal control.

\section{Problem setting} \label{sctn:setting}
\begin{figure}[htb]
\centering
\includegraphics[width=0.35\linewidth]{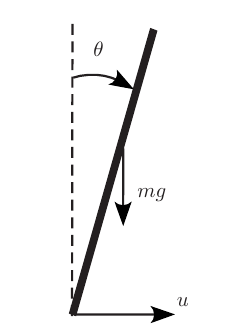}
\caption{2-dimensional pendulum of length $2l$}
\label{fig:pnd_2dim}
\end{figure}
Let us consider a control system for a pendulum with length $l$ and mass $m$ on a mass-less cart as in Fig.~\ref{fig:pnd_2dim}. 
The equation of motion for the pendulum is 
\[
\ddot{\theta}
=\frac{mg\sin\theta-ml\dot{\theta}^2\sin\theta\cos\theta-\cos\theta\,u}{ml/3+ml\sin^2\theta},
\]
where translational motion is ignored. 
With $\bm{x}=\bmat{x_1&x_2}^\top=\bmat{\theta&\dot\theta}^\top$, the control system is written as 
\begin{equation}
\dot{\bm{x}}=f(\bm{x})+g(\bm{x})u \label{eqn:control_system}
\end{equation}
with 
\begin{align*}    
f(\bm{x})&=\bmat{
x_{2}\\[1ex] 
\dfrac{
g\,\sin\left(x_{1}\right)-\frac{l}{2}\,{x_{2}}^2\,\sin\left(2\,x_{1}\right)}{l\,\left({\sin\left(x_{1}\right)}^2+\frac{1}{3}\right)} 
},\\
g(\bm{x})&=\bmat{
0\\[1ex] 
\dfrac{-\cos\left(x_{1}\right)}{ml\left({\sin^2\left(x_{1}\right)}+\frac{1}{3}\right)} 
}.
\end{align*}
Note that $f(0,0)=f(\pi,0)=0$. 
We consider the following infinite-horizon optimal control problem, where 
\begin{equation}
J=\int_0^\infty u(t)^2\,dt \label{eqn:cost}
\end{equation}
is to be minimized 
for (\ref{eqn:control_system}) under the condition that the control trajectory satisfies $\lim_{t\to\infty}\bm{x}(t)=0$ ({\em stable optimal regulator problem}), which enforces the convergence of the pendulum to the upright position. The problem can be interpreted as a minimum-energy swing-up problem under an asymptotic terminal constraint.

A sufficient condition for this problem is characterized by a stabilizing solution $V(x)$ with $\bm{p}^\top =\frac{\partial V}{\partial \bm{x}}=[p_1\ p_2]$ for a partial differential equation (\gls{hjbe})
\begin{equation}
\tilde{H}(\bm{x},\bm{p})=\bm{p}^\top f(\bm{x})-\frac{1}{4}\bm{p}^\top g(\bm{x})g(\bm{x})^\top \bm{p}=0.\label{eqn:hje}
\end{equation}
An approach for an approximate solution for (\ref{eqn:hje}) is to derive an associated Hamiltonian system 
\begin{equation}
    \frac{d\bm{x}}{dt}=\frac{\partial \tilde{H}}{\partial \bm{p}},\ 
    \frac{d\bm{p}}{dt}=-\frac{\partial \tilde{H}}{\partial\bm{x}}\label{eqn:Ham_canonicaleq}
\end{equation}
and compute the stable manifold of the hyperbolic equilibrium (see, e.g., \cite{Lukes:69:sicon,vanderSchaft:92:ac,Sakamoto:08:ieee-ac}). The equation (\ref{eqn:Ham_canonicaleq}) also appears as a necessary condition for the optimal control. 

 We find that the linearization of (\ref{eqn:Ham_canonicaleq}) at $(\bm{x},0)$ is given as below
\[
\bmat{Df(\bm{x})& -\frac{1}{2}g(\bm{x})g(\bm{x})^\top\\
0& -Df(\bm{x})^\top
}.
\]
Here, in our specific problem, due to the special form of the cost functional $J$ and $f(0,0)=f(\pi,0)=0$, the Hamiltonian system has two equilibrium points at $(\bm{x},\bm{p})=(0,0,0,0), (\pi,0,0,0)$. Noting that 
\[
Df(0,0)=\bmat{0&1\\ 3g/l&0},\ Df(\pi,0)=\bmat{0&1\\ -3g/l&0}, 
\]
the linearizations of (\ref{eqn:Ham_canonicaleq}) have unique distributions of eigenvalues as in Table~\ref{tbl:eigenvalues}, and we summarize their importance below: 
\begin{remark}
\begin{enumerate}
\item 
The equilibrium $(0,0,0,0)$ is hyperbolic with 2 stable and 2 unstable eigenvalues, and a 2-dimensional stable manifold exists. The optimally stabilized solutions move along this manifold. On the other hand, all eigenvalues of the equilibrium $(\pi,0,0,0)$ are on the imaginary axis, and the linearization argument cannot be employed for its stability analysis. 
\item
The imaginary eigenvalues appear duplicated, which prevents the direct application of Lyapunov's center theorem\cite[Theorem~5.6.7]{Abraham:79:FOM} for the existence of periodic orbits. This case is notorious for its difficulty in the stability analysis, as even arbitrarily small perturbations or nonlinear effects can unfold the degeneracy and produce instability. It is called a 1:1 or 1:$-1$ resonant equilibrium, and its analysis is a central topic in celestial mechanics, such as Halo orbit bifurcations (see, e.g., \cite{gomez2001dynamics3,Hansman2019jde}). 
\item 
It is particularly interesting that this type of equilibrium appears generically in our optimal control problem, whereas in celestial mechanics it appears with certain parameter values as a bifurcation point. 
\end{enumerate}
\end{remark}
%
\begin{table}[htb]
    \begin{tabular}{c|c|c}
        Equilibrium & $(\bm{x},\bm{p})=(0,0,0,0)$ & $(\bm{x},\bm{p})=(\pi, 0,0,0)$\\[1ex]\hline
        \rule{0pt}{2.5ex}
         Eigenvalues & $\pm\sqrt{3g/l}$, $\pm\sqrt{3g/l}$ & $\pm j\sqrt{3g/l}$, $\pm j\sqrt{3g/l}$\\[1ex]
         & (hyperbolic) & (elliptic)
    \end{tabular}
    \caption{Eigenvalues of the Hamiltonian matrix at two equilibrium points}
    \label{tbl:eigenvalues}
\end{table}
\vskip -1ex
\noindent{\em Problem statement:}\\
In this paper, we aim to characterize the optimally controlled behavior of the pendulum, especially from the hanging position $(\pi, 0,0,0)$. However, it is an equilibrium point of the associated Hamiltonian system for which even its stability is not known. If the optimal control exists, the controlled motion is a connecting trajectory from $(\bm{x},\bm{p})=(\pi, 0,0,0)$ to $(\bm{x},\bm{p})=(0, 0,0,0)$, which is called a heteroclinic orbit. We also wish to understand under what dynamical system properties such an orbit can occur. Summarizing, 
\begin{problem}\label{pblm:stability}
Determine the stability of the Hamiltonian system (\ref{eqn:Ham_canonicaleq}) at $(\bm{x},\bm{p})=(\pi,0,0,0)$.
\end{problem}
\begin{problem}\label{pblm:connection}
Analyze the connecting orbits in (\ref{eqn:Ham_canonicaleq}) from $(\bm{x},\bm{p})=(\pi, 0,0,0)$ to $(\bm{x},\bm{p})=(0, 0,0,0)$. 
\end{problem}
\vskip 1ex
The above problems are directly related to the optimal control objective. If an optimal control exists, the corresponding trajectory must be a solution of the Hamiltonian system (\ref{eqn:Ham_canonicaleq}) that connects the equilibrium $(\pi,0,0,0)$ to $(0,0,0,0)$. However, since the linearization at $(\pi,0,0,0)$ is completely degenerate, consisting of purely imaginary eigenvalues with multiplicity, the standard linearization-based arguments are inconclusive. This necessitates a deeper nonlinear analysis of the Hamiltonian dynamics. As will be seen, the dynamics near this equilibrium exhibit a highly degenerate structure, leading to non-hyperbolic instability mechanisms that differ fundamentally from standard cases.

\section{Numerical simulations}\label{sctn:simulation}
In this section, we present numerical simulations that show a striking feature of the control problem (\ref{eqn:control_system})-(\ref{eqn:cost}) based on \cite{horibe-masterE}. We observe that there exist a number of connecting orbits in (\ref{eqn:Ham_canonicaleq}) from the hanging position to the upright position using trajectories with varying numbers of oscillations, each associated with different control costs. 
 Let $J_k$ denote the cost corresponding to a trajectory with $k$-swings. Representative values are: 
 \begin{gather*}
   J_1=7.1\times 10^{-3}, \  J_3=1.9\times10^{-3}, \    J_5=1.1\times10^{-3}, \\      J_{11}=5.0\times10^{-4},\ \cdots, \ J_{107}=4.7\times10^{-5}. 
 \end{gather*}
 These results indicate that the control cost decreases monotonically as the number of swings increases and suggest that long-time trajectories can achieve the required state transfer with arbitrarily small control effort. 

\begin{figure}[htb]
\centering
\includegraphics[width=1.0\linewidth]{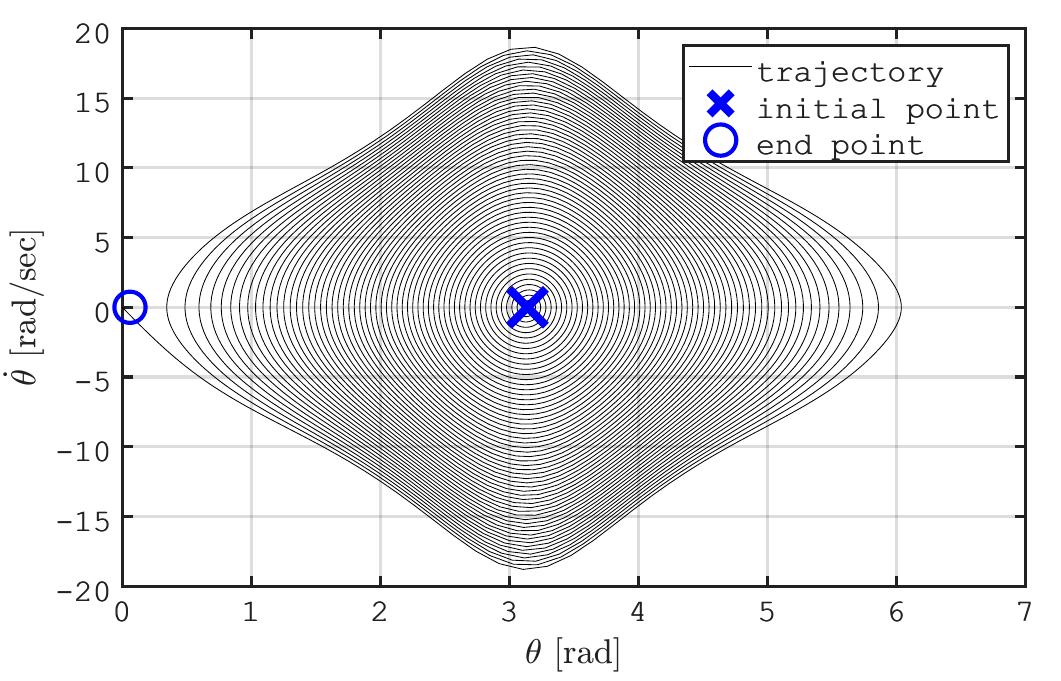}
\caption{Phase portrait of a controlled trajectory exhibiting 107 oscillations before reaching the upright equilibrium. Such long-time trajectories correspond to significantly reduced control cost.}
\label{fig:107swings}
\end{figure}

 Figure~\ref{fig:107swings} illustrates the representative trajectory with 107 swings, showing gradual evolution toward the upright position. In Section, we show that there exists a family of periodic orbits in the phase space, and this trajectory evolves by transitioning from one periodic orbit to another along this family. 
These observations motivate a detailed analysis of the stability properties of the Hamiltonian system (Section~\ref{sctn:stabilityAnal}, Problem~\ref{pblm:stability}) and a variational analysis of the family of the periodic orbits (Section~\ref{sctn:connection}, Problem~\ref{pblm:connection}). 
\section{Stability analysis}\label{sctn:stabilityAnal}
To analyze the stability of $(\bm{x},\bm{p})=(\pi,0,0,0)$, we introduce the change variables $\bm{y}=\bmat{x_1-\pi & x_2}^\top $, obtaining the new Hamiltonian (left side of (\ref{eqn:hje}))
\begin{multline}
H(\bm{y},\bm{p})=p_{1}\,y_{2}
-\dfrac{%
p_{2}\,\left(g\,\sin\left(y_{1}\right)
+\frac{l}{2} {y_{2}}^2\,\sin\left(2\,y_{1}\right)\right)}{l\,\left({\sin^2\left(y_{1}\right)}+\frac{1}{3}\right)}\\
-\dfrac{%
{p_{2}}^2{\cos^2\left(y_{1}\right)}}
{4\,l^2\,m^2\,{\left({\sin^2\left(y_{1}\right)}+\frac{1}{3}\right)}^2}.\label{eqn:orginal_hamiltonian_yp}
\end{multline}
As indicated by the eigenvalues, our analysis starts with the linearized Hamiltonian, and for simplicity of notation, we use a simple parameter normalization to get a Hamiltonian 
\begin{multline}
H(\bm{y},\bm{p})= p_{1}\,y_{2}
-\frac{p_{2}\,\left(\frac{\sin\left(y_{1}\right)}{3}+\frac{{y_{2}}^2\,\sin\left(2\,y_{1}\right)}{2}\right)}{{\sin\left(y_{1}\right)}^2+\frac{1}{3}}\\
-\frac{{p_{2}}^2\,{\cos\left(y_{1}\right)}^2}{18\,{\left({\sin\left(y_{1}\right)}^2+\frac{1}{3}\right)}^2}\label{eqn:ham_normalized}
\end{multline}
and its Hamiltonian matrix
\begin{equation}
A:=J_{sp}\nabla H =\bmat{0 & 1 & 0 & 0\\ -1& 0 & 0& -1\\
0& 0& 0& 1\\ 0&0&-1 &0}\label{eqn:HamMatrix}
\end{equation}
with eigenvalues $\pm j, \pm j$, where $J_{sp}=\left[\begin{smallmatrix}
    0 & I_2\\ -I_2& 0
\end{smallmatrix}\right]$. 
\subsection{Transformation of the linear part}\label{sbsctn:linearTrans}
First, we observe that all the second-order terms in (\ref{eqn:ham_normalized}) corresponding to the linear Hamiltonian matrix (\ref{eqn:HamMatrix}) are given by 
\[
H_0(\bm{x},\bm{p})=x_2p_1-x_1p_2+\frac{1}{2}p_2^2.
\]
In applying Theorem~\ref{thm:app_expansion} in Appendix~\ref{appdx:hamperturb}, the first step is to find conserved quantities along the transpose of the linearized system of (\ref{eqn:HamMatrix}). However, since it is not in the standard form of a Hamiltonian matrix (see \cite[\S~5.4]{Meyer2017i}), finding the conserved quantities for the transposed dynamics is not straightforward. It is therefore convenient to transform the matrix into the standard form so that conserved quantities are well-known. The procedure for finding a symplectic matrix to transform a non-standard Hamiltonian matrix into a standard one is developed in \cite{Laub:74:celes_mech}. It is rather involved, especially for a Hamiltonian matrix whose eigenvalues are all zero or consist of a single pair on the imaginary axis (degenerate cases). Therefore, we present the resulting transformation directly. One can directly confirm that 
\[
P=\bmat{-\frac{\sqrt{2}}{4}& 0& 0& \frac{\sqrt{2}}{2}\\
0& \frac{\sqrt{2}}{4}& -\frac{\sqrt{2}}{2}& 0\\
0 & -\sqrt{2}& 0& 0\\
\sqrt{2}& 0& 0& 0
}
\]
is a symplectic matrix (namely, $P^\top J_{sp} P=J_{sp}$) and that
\[
P^{-1}AP= 
\bmat{0& 1& 0& 0 \\-1& 0& 0& 0\\ 1& 0& 0& 1\\0& 1& -1& 0
},
\]
which is in the standard form. Let $\bm{z}=P^{-1}\left[\begin{smallmatrix}\bm{x}\\\bm{p}\end{smallmatrix}\right]$ and abusing notation, let $H(\bm{z})$ be the new Hamiltonian in the $\bm{z}$-coordinates. Thanks to the symplectic coordinates and the standard form, we have 
\[
H^0(\bm{z}) =z_2z_3-z_1z_4-\frac{1}{2}(z_1^2+z_2^2)
\]
and conserved quantities 
\[
\Gamma_1=z_2z_3-z_1z_4, \quad\Gamma_3=z_3^2+z_4^2
\]
for the transposed linear dynamics $(P^{-1}AP)^\top$. These functions constitute Sokol'skii's symplectic coordinates and lead to the corresponding normal form (see \cite{Sokol'skii1977,Sokol'skii1978}). 
\subsection{Higher order computations with scaling}\label{sbsctn:higherComp}
We perform a Taylor expansion of of $H(\bm{z})$ around $\bm{z}=0$ up to the fourth order:
\begin{align}
&H(\bm{z})=H^0(\bm{z})+H^2(\bm{z})+\cdots; \nonumber\\
& 
\begin{multlined}
\qquad H^2(\bm{z}) =\frac{19}{12}z_{1}{z_{4}}^3
-\frac{3}{2}z_{1}{z_{3}}^2z_{4}
+\frac{3}{2}z_{1}z_{2}z_{3}z_{4}\\
\hspace{4em}+\frac{3}{4}{z_{1}}^2{z_{3}}^2
-\frac{3}{4}{z_{1}}^2z_{2}z_{3} 
+\frac{3}{16}{z_{1}}^2{z_{2}}^2\\
-\frac{37}{16}{z_{1}}^3z_{4}
+\frac{9}{8}z_1^2z_4^2
+\frac{65}{96}{z_{1}}^4.
\end{multlined}\label{eqn:H4}
\end{align}
Note that the original Hamiltonian $H(\bm{x},\bm{p})$ in (\ref{eqn:orginal_hamiltonian_yp}) contains no odd-order terms. 

We now further apply a {\em symplectic scaling} ({\em blow-up transformation}) with a parameter $\varepsilon$: 
\[
\bmat{z_1\\z_2\\z_3\\z_4}
=\bmat{\varepsilon^2z'_1\\\varepsilon^2z'_2\\\varepsilon z'_3\\\varepsilon z'_4} \ \text{or}\ 
\bm{z}'=\bmat{\varepsilon^{-2}I_2& 0\\0& \varepsilon^{-1}I_2}\bm{z}.
\]
For the purpose of stability analysis, we retain terms up to $\varepsilon^2$-order; 
\begin{align}
&H(\bm{z}) = z_2z_3-z_1z_4 -\frac{\varepsilon}{2}(z_1^2+z_2^2)\nonumber\\
&\hspace{6em}+\varepsilon^2\left(\frac{19}{12}z_1z_4^3-\frac{3}{2}z_1z_3^2z_4\right)+O(\varepsilon^3). \label{eqn:ham_upto_3_z}
\end{align}
\subsection{Computation of normal forms with Lie transformation}\label{sbsctn:normalForm}
We introduce Sokol'skii's symplectic coordinates;
\begin{equation}
    \begin{aligned}
z_3&=r\cos\theta, & R&=\frac{1}{r}(z_1z_3+z_2z_4)\\
z_4&=r\sin\theta, & \Theta&=z_2z_3-z_1z_4.
\end{aligned} \label{eqn:sokolskii_sympCoordi}
\end{equation}
One can verify that $(R,\Theta, r,\theta)$ forms a symplectic coordinate system, since 
\[
dR\wedge dr+d\Theta\wedge d\theta=dz_1\wedge dz_3+dz_2\wedge dz_4. 
\]
Note that the conserved quantities $\Gamma_1$, $\Gamma_3$ are 
\[
\Gamma_1=\Theta, \quad \Gamma_3=r^2. 
\]
From (\ref{eqn:app_conserv_original}) in Theorem~\ref{thm:app_expansion} of Appendix~\ref{appdx:hamperturb}, we search a new $H^2$ that depends only on $\Theta$ and $r$, which is our notion of "simple" (see Appendix~\ref{appdx:hamperturb}). The computation procedure in Appendix~\ref{appdx:hamperturb} is written with Poisson brackets which are canonically invariant. So, we may work in the coordinates $(R,\Theta, r,\theta)$. In this case, the Poisson bracket for functions $u,\,v$ is
\begin{equation}
\{u,v\}=
\frac{\partial u}{\partial R}\frac{\partial v}{\partial r}-\frac{\partial u}{\partial r}\frac{\partial v}{\partial R}
+\frac{\partial u}{\partial \Theta}\frac{\partial v}{\partial \theta}-\frac{\partial u}{\partial \theta}\frac{\partial v}{\partial \Theta}.\label{eqn:poisson_sokolskii}
\end{equation}

The Hamiltonian $H(\bm{z})$ in (\ref{eqn:ham_upto_3_z}) is now rewritten in the new coordinates; 
\[
H=\Theta -\varepsilon\left( R^2+\frac{\Theta^2}{r^2}\right)+\varepsilon^2\tilde{H},
\]
 where
\begin{multline*}
\tilde{H}=
- \frac{13}{32} \Theta r^2+\frac{19}{24} \Theta r^2 \cos(2 \theta)  - \frac{37}{96} \Theta r^2 \cos(4 \theta) \\
+ \frac{1}{48}R r^3 \sin(2 \theta) - \frac{37}{96} R r^3 \sin(4 \theta).
\end{multline*}
Here, we have the following correspondence with the procedure in Appendix~\ref{appdx:hamperturb};
\[
H_0^0=\Theta, \quad H_1^0=-R^2-\frac{\Theta^2}{r^2},\quad
H_0^2=\tilde{H}.
\]
In the first equation (\ref{eqn:app_Lie1}) for $H_0^1$, $H_1^0$ corresponds to the linear dynamics and there is no need to simplify it. Therefore we take $W_1=0$. This allows to combine (\ref{eqn:app_Lie2}) and (\ref{eqn:app_Lie3}) to get 
\[
H_0^2=H_2^0+\{H_0^0,W_2\}.
\]
For the simplification, we note that from (\ref{eqn:poisson_sokolskii})
\[
\{H_0^0,W_2\}=\frac{\partial W_2}{\partial \theta}
\]
and choose $W_2$ so as to eliminate the terms involving $\theta$. The resultant $H_0^2$ is
\[
H_0^2=-\frac{13}{32} \Theta r^2. 
\]
The resulting normal form, up to $\varepsilon^2$-order, is given by
\begin{equation}
H=\Theta -\varepsilon\left( R^2+\frac{\Theta^2}{r^2}\right)-\frac{13}{32}\varepsilon^2\Theta r^2+O(\varepsilon^3).\label{eqn:normalForm} 
\end{equation}
This is a form of Sokol'skii's normal form in the sense that higher order terms in $H$ are written with conserved quantities (for the transposed linear dynamics) $r^2$ and $\Theta$. 
\subsection{Proof of instability}\label{sbsctn:chetaev}
The Hamiltonian system corresponding to (\ref{eqn:normalForm}) is given by
\begin{equation*}
\left\{
\begin{aligned}
\dot R &=\frac{2\varepsilon\Theta^2}{r^3}
    -2\varepsilon^3\Theta r+O(\varepsilon^3)\\
\dot\Theta &=O(\varepsilon^3)\\
\dot r &= 2\varepsilon R+O(\varepsilon^3)\\
\dot\theta &= -1+\frac{2\varepsilon\Theta}{r^2}
    +\alpha\varepsilon^2r^2+O(\varepsilon^3),
\end{aligned}\right. 
\end{equation*}
where $\alpha=\frac{13}{32}$. 
Applying the inverse transformation of (\ref{eqn:sokolskii_sympCoordi}), the system can be expressed in the $\bm{z}$-coordinates as 
\begin{equation}
\left\{
\begin{aligned}
\dot{z}_1 &= z_2-\alpha\varepsilon^2(z_2r^2+2z_3\Theta)+O(\varepsilon^3)\\
\dot{z}_2 &= -z_1
+\alpha\varepsilon^2(z_1r^2 -2z_4\Theta)+O(\varepsilon^3)\\
\dot{z}_3 &= z_4
+2\varepsilon z_1
-\alpha\varepsilon^2z_4r^2+O(\varepsilon^3)\\
\dot{z}_4 &= -z_3
+2\varepsilon z_2
+\alpha\varepsilon^2z_3r^2+O(\varepsilon^3). 
\end{aligned}\right.\label{eqn:hamsys_z_normalform}
\end{equation}
As a candidate Chetayev function, we take \[V=z_1z_3+z_2z_4\] and its time derivative along the trajectories of (\ref{eqn:hamsys_z_normalform}) is 
\[
\dot V =2\varepsilon(z_1^2+z_2^2)-2\alpha\varepsilon^2(z_3^2+z_4^2)\Theta +O(\varepsilon^3).
\]
Define an open set $\Omega\subset \mathbb{R}^4$ as
\begin{multline*}
\Omega=\{\bm{z}\in\mathbb{R}^4\,|\, 
V=z_1z_3+z_2z_4>0, \\
\Theta=z_2z_3-z_1z_4<0\}. 
\end{multline*}
Conditions (i) and (ii) in Theorem~\ref{thm:app_chetaev} in Appendix are satisfied, since $(0,0,0,0)\in \partial\Omega$ and $V>0$ on $\Omega$. Also, Condition~(iii) is $V=0$ on $\partial\Omega$ and is satisfied. For Condition~(iv), since $\Theta<0$ on $\Omega$ and $\alpha=13/32>0$, we have$\dot V>0$ on $\Omega$ for a sufficiently small $\varepsilon>0$. Here, we note that $O(\varepsilon^3)$-order terms in $\dot V$ are polynomials of degree~4, and do not affect the sign of $\dot V$. We thus obtain the following result (see Problem~\ref{pblm:stability}). 
\vskip 1ex
\begin{proposition}\label{prop:instability}
The equilibrium $(\bm{x},\bm{p})=(\pi,0,0,0)$ of Hamiltonian system (\ref{eqn:Ham_canonicaleq}) is unstable. 
\end{proposition}

\begin{figure}[htb]
    \centering
    \includegraphics[width=0.95\linewidth]{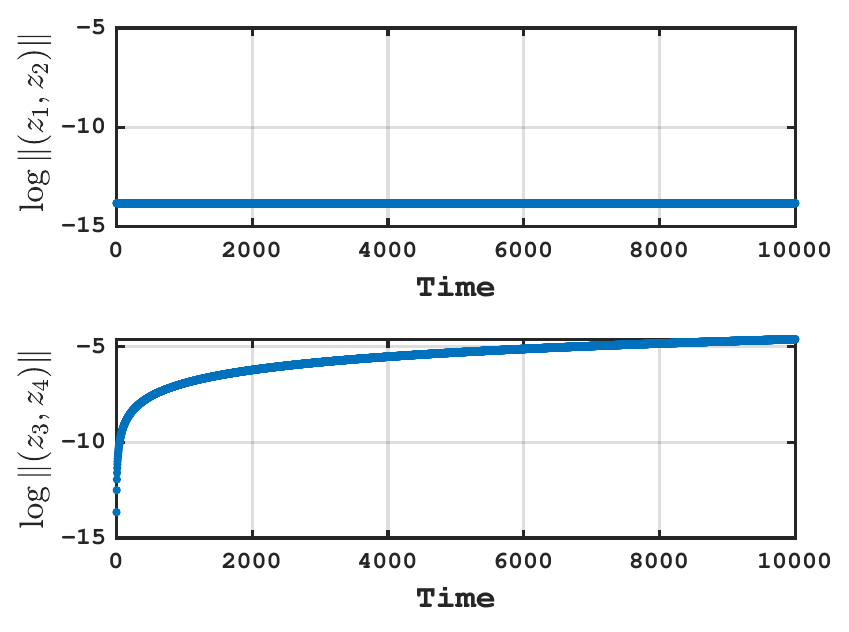}
    \caption{Numerical simulation of system  (\ref{eqn:hamsys_z_normalform}) with initial condition $\bm{z}(0)=(0,10^{-6},10^{-6},0)$. The trajectory exhibits slow divergence from the equilibrium, illustrating the weak (non-hyperbolic) instability predicted by Proposition~\ref{prop:instability}.}
    \label{fig:instability}
\end{figure}

\vskip 1ex
Fig.~\ref{fig:instability} illustrates a numerical simulation of (\ref{eqn:hamsys_z_normalform}) with initial condition $(0,10^{-6},10^{-6},0)$, demonstrating weakly unstable behavior. 
\section{Mechanism of connection between hanging and upright equilibria}\label{sctn:connection}
\subsection{Conserved quantity and periodic orbits} 
The Hamiltonian system (\ref{eqn:Ham_canonicaleq}) in $(\bm{y},\bm{p})$-coordinates is given as
\begin{subequations}
\begin{align}
    &\frac{d\bm{y}}{dt}=f(\bm{y})-\frac{1}{2}\tilde{g}(\bm{y})\tilde{g}(\bm{y})^\top \bm{p}\\
    &\frac{d\bm{p}}{dt}=-\left(\frac{\partial f}{\partial \bm{y}}\right)^\top\bm{p}
    +\frac{1}{4}\frac{\partial}{\partial \bm{y} }\left[\bm{p}^\top \tilde{g}(\bm{y})\tilde{g}(\bm{y})^\top \bm{p}\right],
\end{align} \label{eqn:Hamsys-y}
\end{subequations}
where $\tilde{g}(\bm{y})=\bmat{0 & \dfrac{\cos\left(y_{1}\right)}{\left({\sin^2\left(y_{1}\right)}+\frac{1}{3}\right)}}^\top$. It is readily seen that 
\[
\frac{d\bm{y}}{dt}=f(\bm{y}),\quad \bm{p}=0
\]
is always a solution for (\ref{eqn:Hamsys-y}), and the first equation above has a conserved quantity 
\begin{equation}
E(\bm{y})=\frac{1}{2}(\frac{1}{3}+\sin^2(y_1))y_2^2-\frac{1}{3}\cos(y_1).\label{eqn:conserved}
\end{equation}
This means that (\ref{eqn:Hamsys-y}) has a family of periodic orbits $(\bm{y}(t),0)$ satisfying $E(\bm{y}(t))=E$, where $E$ is a constant determined by the initial condition $\bm{y}(0)$. 

The period of each periodic orbit is obtained by a quadrature as is shown below: 
\begin{equation}
T(A)=2\sqrt{6}\int_0^A \sqrt{
    \frac{1/3+\sin^2(y_1)}{\cos(y_1)-\cos(A)}
    }\, dy_1, \label{eqn:period}
\end{equation}
where $A\in(0,\pi)$ is the amplitude of the periodic orbit. 
It can be seen that $T(A)\to0$ as $A\to0$ and $T(A)\to\infty$ as $A\to\pi$. Thus, we have shown the following. 
\begin{proposition}
In the Hamiltonian system (\ref{eqn:Hamsys-y}), there exists a family of periodic orbits with $\bm{p}=0$ centered at $(0,0,0,0)$. It is parametrized with the pendulum oscillation amplitude $A\in(0,\pi)$, and as $A\to\pi$, the periodic orbit converges to the homoclinic orbit connecting $(\bm{y},\bm{p})=(0,0,0,0)$ and $(\bm{y},\bm{p})=(\pi,0,0,0)$. 
\end{proposition}
\vskip 1ex
In other words, we have obtained an invariant cylinder $\Gamma$ as
\[
\Gamma=\left\{
(y_1,y_2,0,0)\,|\, E(\bm{y})=-\frac{1}{3}\cos(A),\ A\in(0,\pi)
\right\}.
\]
The proposition is consistent with the fact that $\bm{p}=0$ corresponds to the uncontrolled dynamics, and the periodic orbits are inherent oscillations of the pendulum itself. 
\subsection{Variational analysis of the periodic orbits} The goal of this section is to understand how control inputs can drive the pendulum 
from the hanging position $(\bm{y},\bm{p})=(0,0,0,0)$ to the upright position $(\bm{y},\bm{p})=(\pi,0,0,0)$ through the dynamics of (\ref{eqn:Hamsys-y}). In Proposition~\ref{prop:instability}, we have shown that one boundary of $\Gamma$, which is formed by periodic orbits of the uncontrolled system, is unstable. This suggests that even small perturbations in the costates $\bm{p}$ (or equivalently, small control inputs) may induce trajectories that depart from this equilibrium and evolve in the vicinity of $\Gamma$. 
This observation raises the possibility of motions that start from the hanging position, evolve along the family of periodic orbits, and eventually reach the upright equilibrium along its stable manifold (Figure~\ref{fig:107swings}). To explore this mechanism, we investigate the linearization of each periodic orbit in the family. As we will see, the periodic orbits exhibit a highly nonstandard and maximally degenerate Floquet structure, which plays a crucial role in this mechanism. 

Let $z=\bmat{\bm{y}^\top,\bm{p}^\top}^\top$ and rewrite (\ref{eqn:Hamsys-y}) as $\dot{z}=f(z)$. Let $\varphi(z,t)$ be its solution starting from $z$ at $t=0$. Then, the variational equation of (\ref{eqn:Hamsys-y}) is a $4\times 4$-matrix-valued equation
\begin{equation}
\dot{\Phi}(t)=Df(\varphi(t,z))\Phi(t). \label{eqn:variationaleq}
\end{equation}
The scheme we employ for this variational analysis is the 2-stage Gauss-Legendre collocation method (order 4) (see, e.g., \cite{hairer:06:GNISPAODE}), which is a symplectic integrator suitable for the analysis of Hamiltonian systems. 

\subsubsection{Computation of period $T(A)$}
\begin{figure}[htb]
\centering
\includegraphics[width=0.9\linewidth]{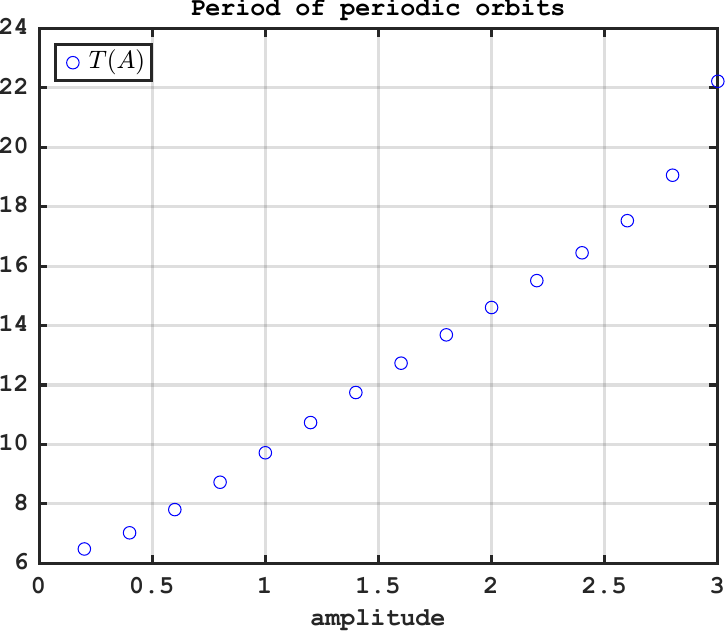}
\caption{Period $T(A)$ of the periodic orbits forming the invariant cylinder $\Gamma$, plotted against the oscillation amplitude $A$. }
\label{fig:period}
\end{figure}

In the formula for $T(A)$ in (\ref{eqn:period}), there are two singularities, one of which can be removed by the change of variables
\[
\sin(y_1/2)=k\sin \phi,\ k=\sin(A/2),\ \phi\in[0,\pi/2]. 
\]
Thus, we compute $T(A)$ by
\[
T(A)=4\sqrt{3}\bigint_0^{\dfrac{\pi}{2}}\!\!\!\frac{\sqrt{1/3+4k^2\sin^2\phi(1-k^2\sin^2\phi)}}{\sqrt{1-k^2\sin^2\phi}}\,d\phi.
\]
Figure~\ref{fig:period} shows that the period $T(A)$ monotonically increases as the amplitude $A$ increases. 
\subsubsection{Monodromy matrices of periodic orbits}
The monodromy matrix $M(A)\in\mathbb{R}^{4\times4}$ of a periodic orbit in (\ref{eqn:Hamsys-y}) starting at $(A,0,0,0)$ is the solution of (\ref{eqn:variationaleq}) at $t=T(A)$ with $z=(A,0,0,0))$ and $\Phi(0)=I_4$. Actual computation for the monodromy matrix is to solve a system of $4+4^2$ equations consisting of (\ref{eqn:Hamsys-y}) and (\ref{eqn:variationaleq}) with the initial condition $(A,0,0,0, I_4)$. 

\begin{figure}[htb]
\centering
\includegraphics[width=0.9\linewidth]{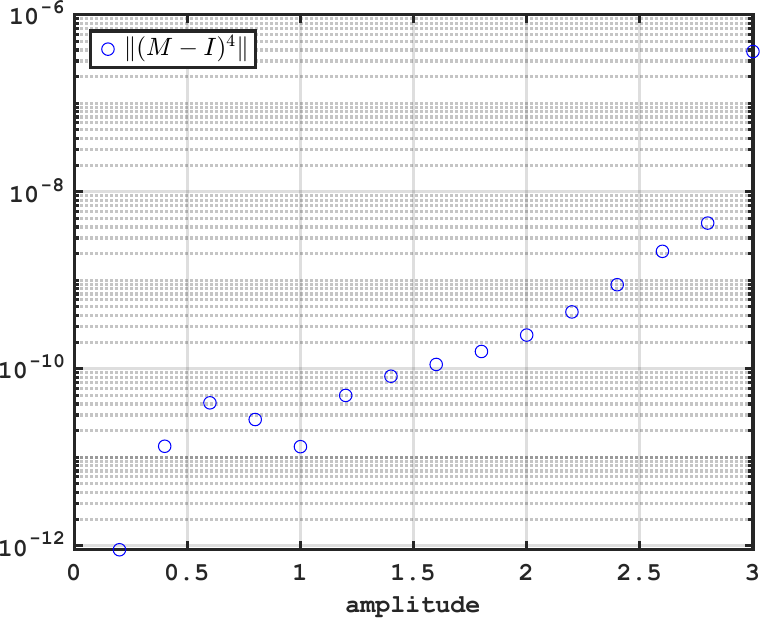}%
\caption{Norm of $(M(A)-I_4)^4$ versus amplitude $A$. The near-vanishing values confirm the approximate nilpotency of order four.}
\label{fig:eig_1}
\end{figure}

Our computations reveal that the monodromy matrix exhibits a maximally degenerate structure, summarized as follows. 
\begin{enumerate}[(i)]
    \item All four Floquet multipliers of $M(A)$ are equal to 1. 
    \item The invariant cylinder $\Gamma$ is not normally hyperbolic, as there is no exponential splitting in the transversal directions.
    \item The eigenvalue 1 is associated with a single Jordan block of size 4 (maximal length).  
    \item This degeneracy persists across the entire family of periodic orbits, i.e., for all amplitude $A$.  
\end{enumerate}
\vskip 1ex
Item~(i) is seen in Figure~\ref{fig:eig_1} which depicts $\|(M(A)-I_4)^4\|$ against $A$. 
Item~(ii), in contrast to the standard situation, indicates the absence of exponentially attracting and repelling dynamics in the vicinity of $\Gamma$. Note that if there is a periodic orbit in a Hamiltonian system, the corresponding Monodromy matrix has at least two Floquet multipliers equal to 1 (one tangent direction of the orbit, the other, tangent direction for the family). See, for instance, \cite{BeckSakamoto2026} for periodic orbits of Hamiltonian systems arising in an optimal control problem that exhibit a hyperbolic structure. 

\begin{figure}[htb]
\centering
\includegraphics[width=0.9\linewidth]{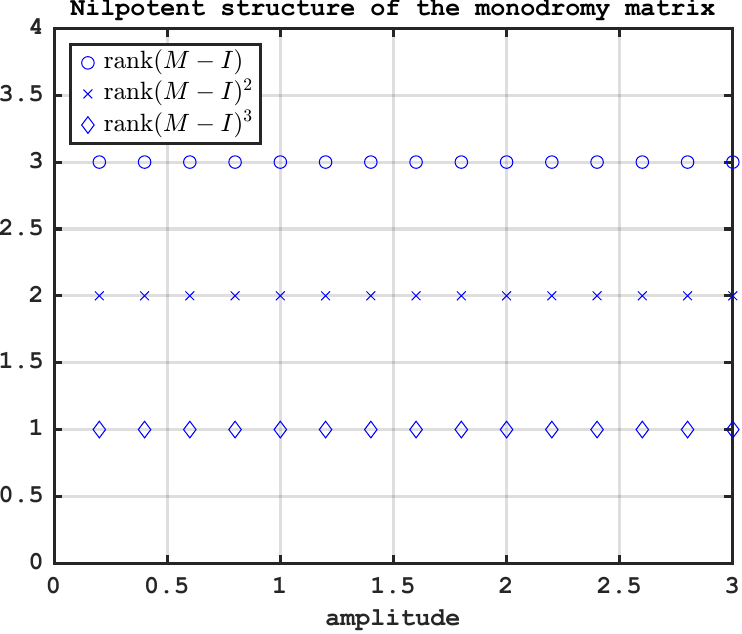}%
\caption{Rank of successive powers of $M(A)-I_4$ along the family of periodic orbits. The cascade confirms that the eigenvalue 1 is associated with a single Jordan block of size 4. }
\label{fig:nilpotent}
\end{figure}

Item~(iii) is confirmed by $(M(A)-I_4)^4=0$ while $(M(A)-I_4)^3\neq0$, which is depicted in Figure~\ref{fig:nilpotent}. This is highly nonstandard: in autonomous Hamiltonian systems, periodic orbits typically have a semisimple multiplier 1 of multiplicity two, and higher-order degeneracies of this type are not observed in generic settings. 
This maximal degeneracy has a direct dynamical consequence. Writing $N=M(A)-I_4$, the relation $N^4=0$ implies that iterates of the monodromy matrix grow polynomially in time, with terms up to third-order. Thus, perturbations evolve via slow, non-exponential drift along the invariant cylinder $\Gamma$. This behavior is consistent with the weak instability observed in \S~\ref{sctn:stabilityAnal} and the long-time trajectories in \S~\ref{sctn:simulation}. 

\begin{figure}[htb]
\centering
\includegraphics[width=0.9\linewidth]{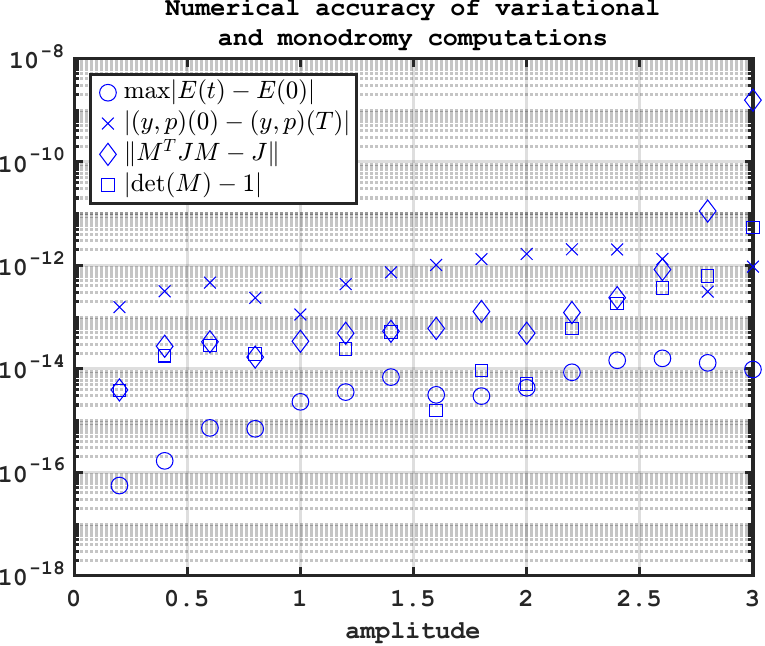}%
\caption{The numerical accuracy of the variational and monodromy computations. The deviations in the conservation law, periodicity, symplecticity, and determinant.  }
\label{fig:accuracy}
\end{figure}

Numerical accuracy is verified via conservation of energy, periodicity, symplecticity, and determinant conditions, all satisfied within $10^{-12}$ (Figure~\ref{fig:accuracy}). 
Figure~\ref{fig:accuracy} shows that the monodromy computations are reliable over the range of amplitudes shown in Figures~\ref{fig:eig_1} and \ref{fig:nilpotent}. Near the homoclinic limit $A\to\pi$, the period becomes large, and the variational integration becomes increasingly ill-conditioned, so we do not make claims there based solely on numerics. 
\subsubsection{Stable manifold of the upright equilibrium}
A simple computation shows that the eigenvectors corresponding to the eigenvalue $-1$ (duplicated) at the hyperbolic (upright) equilibrium of (\ref{eqn:Hamsys-y}) are 
\[
\bmat{1&-1&0&0}^\top,\quad \bmat{1&0&2&2}^\top, 
\]
and the stable manifold of the upright equilibrium is tangent to the subspace spanned by them. It is seen that a controlled trajectory approaching the equilibrium along the first eigenvector requires zero control input (the last two elements correspond to $\bm{p}$). This computation is consistent with the existence of small connecting control driving the pendulum from the hanging position to the upright position in \S~\ref{sctn:simulation}.

\section{Conclusions}
This paper analyzed an infinite-horizon optimal control problem for a pendulum via its associated Hamiltonian system. Due to degeneracy of the linearization at the hanging equilibrium, standard linear methods are inconclusive. Using Lie series and Chetaev’s instability theorem, we established weak, non-hyperbolic instability.
We further identified a family of periodic orbits forming an invariant cylinder and showed that the associated monodromy matrices exhibit a maximally degenerate Floquet structure with a single Jordan block at the unit multiplier. This induces slow, non-exponential (polynomial-type) dynamics and provides a mechanism for long-time transitions with arbitrarily small control effort.

Although demonstrated for a pendulum, this phenomenon arises from degeneracy in the optimal control formulation and may occur more broadly. In particular, it suggests that optimal control may fail to exist despite stabilizability. Future work includes establishing this nonexistence rigorously and extending the analysis to more general systems.

{\bf Acknowledgment.} The author would like to thank Fabian Beck for helpful discussions. 
\bibliographystyle{IEEEtran}
%

\appendices
\setcounter{equation}{0}
\renewcommand{\thetheorem}{\thesubsection\arabic{theorem}}
\renewcommand{\theequation}{\thesubsection\arabic{equation}}
\section{Hamiltonian normal forms and Lie series}\label{appdx:hamperturb}
Let $H(\bm{x})$, $\bm{x}\in\mathbb{R}^{2n}$ be analytic and suppose that it is expanded with homogenous polynomials as
\[
H(\bm{x})=H_0(\bm{x})+H_1(\bm{x})+\cdots
    =\sum_{i=0}^\infty H_i(\bm{x}), 
\]
where $H_i(\bm{x})$ is a polynomial of degree $i+2$ and let $A=
J_{sp}\nabla H_0$, where $J_{sp}\left[\begin{smallmatrix}
    0 & I_n\\ -I_n& 0
\end{smallmatrix}\right]$. 
\vskip 1ex
\begin{theorem}\label{thm:app_expansion}
There exists an analytic symplectic transformation $\bm{x}=X(\bm{y})$ with $DX(0)=I_{2n}$ such that $H^\ast (\bm{y}):=H(X(\bm{y}))$ satisfies
\begin{equation}
    H^i(e^{tA^\top}\bm{y})=H^i(\bm{y}) \label{eqn:app_conserv_original}
\end{equation}
for all $i=0,1,\ldots$, $t\in\mathbb{R}$ and $\bm{y}\in\mathbb{R}^{2n}$, where $H^i(\bm{y})$ is a polynomial of degree $i+2$ such that
\begin{equation}
H^\ast(\bm{y})=H^0(\bm{y})+H^1(\bm{y})+\cdots
    =\sum_{i=0}^\infty H^i(\bm{y}).\label{eqn:app_H_series_new}
\end{equation}
\end{theorem}

\vskip 1ex
When the new Hamiltonian $H^\ast (\bm{y})$ is simple in some sense, (\ref{eqn:app_H_series_new}) is said to be in normal form. Being simple may have different meanings depending on the purpose of the analysis, but, the condition (\ref{eqn:app_conserv_original}) plays a key role. 

A systematic way to find a symplectic transformation $X(\bm{y})$ satisfying (\ref{eqn:app_conserv_original}) is to use the Lie series expansion. 
It is assumed that $X(\bm{y})$ is generated by a Hamiltonian system
\[
\frac{d\bm{x}}{d\varepsilon}=J_{sp}\nabla W(\varepsilon,\bm{x}),\quad \bm{x}(0)=\bm{y}.
\]
Further assuming that the original Hamiltonian and $W$ have the following expansions
\begin{align*}
H(\varepsilon,\bm{x})&=\sum_{i=0}^\infty \left( \frac{\varepsilon}{i!}\right)H_i^0(x),\\
W(\varepsilon,\bm{x})&=\sum_{i=0}^\infty \left( \frac{\varepsilon}{i!}\right)W_{i+1}(x),
\end{align*}
$W_1,\, W_2,\,\ldots$ are computed so that each term in the new Hamiltonian 
\[
H^\ast(\varepsilon,\bm{y})=H(\varepsilon,X(\bm{y}))=\sum_{i=0}^\infty \left( \frac{\varepsilon}{i!}\right)H_0^i(y)
\]
is simple. This computation is carried out on a term-by-term basis 
in the recursive computation using Poisson brackets;
\begin{subequations}
\begin{align}
H_0^1 &= H_1^0 +\{H_0^0,W_1\}, \label{eqn:app_Lie1}\\
H_1^1 &= H_2^0 +\{H_1^0,W_1\}+\{H_0^0,W_2\},\label{eqn:app_Lie2}\\
H_0^2 &= H_1^1+\{H_0^1,W_1\}, \label{eqn:app_Lie3}\\
&\cdots.\nonumber
\end{align}
\end{subequations}
\section{Chetayev's instability theorem}
\setcounter{equation}{0}
Let $U\subset\mathbb{R}^n$ be an open set containing the origin and $f:U\to\mathbb{R}^n$ be a smooth vector field with $f(0)=0$. Consider a differential equation 
\begin{equation}
\dot{\bm{z}}=f(\bm{z}). \label{eqn:app_chetaev}
\end{equation}
\begin{theorem}\label{thm:app_chetaev}
Suppose that there exist an open set $\Omega\subset U$ 
 and a $C^1$ function $V:U\to\mathbb{R}$ such that
\begin{enumerate}[(i)]
    \item $0\in \partial\Omega$, where $\partial\Omega$ is the boundary of $\Omega$,
    \item $V(z)>0$ for $z\in\Omega$,
    \item $V(z)=0$ for $z\in\partial\Omega$,
    \item $\dot V:=(\nabla V\cdot f)(z)>0$ for $z\in\Omega$. 
\end{enumerate}
Then, the origin is an unstable equilibrium of (\ref{eqn:app_chetaev}).
\end{theorem}

\end{document}